\documentclass[11pt]{article}
\usepackage{amsmath,amsfonts,amssymb}
\usepackage[utf8]{inputenc}
\usepackage{tikz,ifthen}

\newtheorem{theorem}{Theorem}[section]
\newtheorem{lemma}[theorem]{Lemma}

\newtheorem{corollary}[theorem]{Corollary}

\newcommand{\proof}{\noindent{\bf Proof.\ }}
\newcommand{\qed}{\hfill $\square$ \bigskip}

\begin{document}

\title{Domination game and minimal edge cuts}

\author{
Sandi Klav\v zar $^{a,b,c}$
\and
Douglas F. Rall $^{d}$
}

\date{}

\maketitle

\begin{center}
$^a$ Faculty of Mathematics and Physics, University of Ljubljana, Slovenia\\
$^b$ Faculty of Natural Sciences and Mathematics, University of Maribor, Slovenia\\
$^c$ Institute of Mathematics, Physics and Mechanics, Ljubljana, Slovenia\\
{\tt sandi.klavzar@fmf.uni-lj.si}
\medskip

$^d$ Department of Mathematics, Furman University, Greenville, SC, USA\\
{\tt doug.rall@furman.edu}
\end{center}

\begin{abstract}
In this paper a relationship is established between the domination game and minimal edge cuts. It is proved that the game domination number of a connected graph can be bounded
above in terms of the size of  minimal edge cuts.    In particular, if $C$ a minimum edge cut of a connected graph $G$, then $\gamma_g(G) \le \gamma_g(G\setminus C) + 2\kappa'(G)$.
Double-Staller graphs are introduced in order to show that this upper bound can be attained for graphs with a bridge.  The obtained results are used to extend the family of
known traceable graphs whose game domination numbers are at most one-half their order.   Along the way  two technical lemmas, which seem to be generally applicable for the study of the domination game, are proved.
\end{abstract}


\medskip\noindent
{\bf AMS Subj.\ Class:} 05C57, 05C69, 91A43


\section{Introduction and preliminaries}
\label{sec:intro}

The investigation of the domination game that was initiated in~\cite{BKR2010} generated investigations by a number of authors;~\cite{BDKK-2017, bujtas-2015b, hl-2017,  NSS2016, xu-2019,  XLK-2018} is just a small selection of the papers. Offspring of the game include the total domination game~\cite{henning-klavzar-rall-2015} (see also~\cite{bresar-2017, bujtas-henning-tuza-2016, combinatorica-2017, henning-rall-2016}), the disjoint domination game~\cite{bujtas-tuza-2016}, the transversal game on hypergraphs~\cite{bujtas-henning-tuza-2016, bujtas-henning-tuza-2017, bujtas-2018+}, and the game distinguishing number~\cite{gravier-2017}.

The domination game is played on a graph $G$ by Dominator and Staller. They take turns choosing a vertex from $G$ such that at least one previously undominated vertex becomes dominated. The game is over when $G$ becomes dominated. Dominator wants to minimize the number of vertices played, and Staller wishes to maximize it. A game is a {\em D-game} (resp.\ {\em S-game}) if Dominator (resp.\ Staller) has the first move. If both players play optimally, then the number of moves played is the {\em game domination number} $\gamma_g(G)$ of $G$. The corresponding invariant for the S-game is denoted $\gamma_g'(G)$.

In this paper we are interested in the interplay between the domination game and minimal edge cuts.  In the rest of this section we give the necessary definitions and
recall several known results that we need.   In Section~\ref{sec:two-lemmas} we prove two technical lemmas, both of which provide more insight into the game from Staller's point of view.
Our main result, Theorem~\ref{thm:main1}, which relates the game domination number to the size of  minimal edge cuts, is proved in Section~\ref{sec:edge-cuts}.   An immediate corollary of the theorem is that if $C$ a minimum edge cut of a connected graph $G$, then $\gamma_g(G) \le \gamma_g(G\setminus C) + 2\kappa'(G)$, where $\kappa'(G)$ denotes the edge-connectivity of $G$.  In the subsequent section we prove that the upper bound
in this corollary can be attained for graphs with a bridge.  For this sake the concept of a double-Staller graph is introduced.  In the final section we give an application of the results from Section~\ref{sec:edge-cuts} and use it to extend the family of
known traceable graphs whose game domination numbers are at most one-half their order.

If $x$ is a vertex of a graph $G$, then $N[x] = \{y:\ xy\in E(G)\}\cup \{x\}$ is the {\em closed neighborhood} of $x$. The vertex $x$ {\em dominates} every vertex in $N[x]$. The order of a graph $G$ is denoted by $n(G)$, and for a positive integer $m$ we use the notation $[m] = \{1,\ldots, m\}$.  A graph is called {\em traceable} if it contains a hamiltonian path.

The sequence of moves in a D-game will be denoted with $d_1, s_1, d_2, s_2, \ldots$, and the sequence of moves is an S-game with $s_1', d_1', s_2', d_2', \ldots$ A {\em partially dominated graph} is a graph together with a declaration that some vertices are already dominated, meaning that they need not be dominated in the rest of the game. If $S\subseteq V(G)$, then $G|S$ denotes the partially dominated graph in which vertices from $S$ are already dominated. If $S=\{x\}$, we will abbreviate $G|\{x\}$ to $G|x$. If $G|S$ is a partially dominated graph, then $\gamma_g(G|S)$ and $\gamma_g'(G|S)$ denote the optimal number of moves in the D-Game and the S-Game, respectively, played on $G|S$.  The {\em residual graph} of $G|S$ is the graph obtained from $G$ by removing all vertices $u$ for which each vertex in $N[u]$ is dominated and all edges joining dominated vertices.

The following result is one of the most useful tools for the domination game.

\begin{lemma}
\label{lem:Con-P}{\em (Continuation Principle, \cite{kinnersley-2013})}
Let $G$ be a graph with $A,B\subseteq V(G)$. If $B\subseteq A$, then $\gamma_g(G|A)\leq \gamma_g(G|B)$ and $\gamma^{\prime}_g(G|A)\leq \gamma^{\prime}_g(G|B)$.
\end{lemma}

Using the Continuation Principle, a short proof can be given for the following result from~\cite{BKR2010, kinnersley-2013}.

\begin{theorem}
\label{thm:diff-at-most-1}
If $G$ is a partially dominated graph, then $|\gamma_g(G) - \gamma_g'(G)|\le 1$.
\end{theorem}

\section{Two lemmas}
\label{sec:two-lemmas}

In this section we prove two lemmas to be used later in the paper and that could be of independent interest.

\begin{lemma} \label{lem:prevertex}
If  $G$ is a graph, $S \subseteq V(G)$, and $x \in V(G)$, then $$\gamma_g'(G|S) \leq \gamma_g'(G|(S\cup N[x])) +2\,.$$
\end{lemma}

\proof
Consider the S-game played on $G|S$.  Let $s_1'$ be an optimal first move of Staller, and let $d_1'=x$ provided this is a legal move.  Otherwise, let $d_1'$ be an arbitrary legal move.
From now on Dominator is imagining another S-game played on $G|(S\cup N[x])$.    In the real game played on $G|S$ Staller is playing optimally, but in the imagined game Dominator is playing optimally.  The strategy of Dominator will be such that after each move the set of vertices dominated in the imagined game is
a subset of the set of vertices dominated in the real  game.  Note that this is certainly the case after the moves $s_1'$ and $d_1'$ were made in the real game.  Whenever it is Dominator's turn he plays optimally in the imagined game and copies that move into the real game if it is a legal move there.  Otherwise Dominator plays an arbitrary legal move in the real game.
  Whenever it is Staller's turn she plays optimally in the real game, and Dominator copies this move into the imagined game.  Because of Dominator's strategy described above, this copied
  move of Staller to the imagined game will always be legal.

Let $r$ be the number of moves
made in the real game after $s_1'$ and $d_1'$ are made, and let $t$ be the number of moves played in the imagined game.  By the strategy of Dominator, we get $r \leq t$.    Since Staller is playing optimally in the real game, $\gamma_g'(G|S) \leq 2+r$.  On the other hand, in the imagined game Dominator is playing optimally, which implies that $\gamma_g'(G|(S\cup N[x])) \geq t$.  Combining these inequalities we see
that
$$\gamma_g'(G|S) \leq 2+r \leq 2+t \leq 2+\gamma_g'(G|(S\cup N[x]))\,,$$
which finishes the proof.
\qed

By an easy application of the Continuation Principle, the following corollary, which we will use later in the paper, is immediate.

\begin{corollary}\label{cor:prevertex}
If  $G$ is a graph, $S \subseteq V(G)$, and $x \in V(G)$, then $$\gamma_g'(G|S) \leq \gamma_g'(G|(S\cup \{x\})) +2\,.$$
\end{corollary}

If $p$ is a non-negative integer, then the {\em Staller $p$-pass game} is a D-game in which Staller is allowed not to play, that is, to {\em pass}, at most $p$ times. Clearly, the Staller $0$-pass game is the usual D-game. The number of moves (the pass moves of Staller are not counted) in the Staller $p$-pass game, when both players play optimally,  is denoted ${\gamma}_g^{St,p}(G)$. Our second lemma reads as follows.

\begin{lemma}
\label{lem:pass-game}
If $G$ is a graph and $p$ a non-negative integer, then ${\gamma}_g^{St,p}(G) \le \gamma_g(G) + p$.
\end{lemma}

\proof
There is nothing to be proved if $p=0$, hence assume in the rest that $p\ge 1$.

We consider two games played at the same time. The first one, addressed as the {\em real game}, is the Staller $p$-pass game played on $G$. The second game, the {\em imagined game}, is a D-game played on $G$ which is imagined (and played) only by Dominator. Dominator plays optimally in the imagined game and copies each of his moves to the real game. (This means that Dominator might not play optimally in the real game.) After each move of Dominator in the real game, Staller either plays a move or passes it. In the first case Dominator copies Staller's move from the real game to the imagined game provided this move is legal. The situation that Staller's move is not legal in the imagined game will be treated later. In the second case (when Staller passes a move in the real game), Dominator imagines an arbitrary legal move of Staller in the imagined game and replies optimally afterwards with his own move (and copies it back to the real game).

Note that by the above described strategy, at each stage of the games the set of vertices being already dominated in the real game is a subset of the vertices being already dominated (at the same stage of the games) in the imagined game. It follows that each (optimal) move of Dominator from the imagined game is indeed a legal move in the real game. And, as said, a move of Staller from the real game might not be legal in the imagined game. In any case, the imagined game ends no later than the real game. Let $s$ be the number of moves played in the imagined game.

Suppose that, when the imagined game ends,  Staller has made $q$ pass moves, $q\le p$, in the real game and let $s^{(1)}, \ldots, s^{(q)}$ be the moves of Staller in the imagined game that were imagined by Dominator as the moves corresponding to her pass moves in the real game. Suppose for the time being that every move of Staller from the real game has been a legal move in the imagined game until the imagined game has been finished. Note that until this point, the number of moves played in the imagined game is precisely $q$ larger than the number of moves in the real game. The number of moves played in the real game until this moment is thus $s - q$. The strategy of Dominator in the remainder of the real game is to play the vertices $s^{(1)}, \ldots, s^{(q)}$ during  the next (at most) $2q$ moves of himself and Staller. Of course, if at some point in the remainder of the real game one of these is not a legal move, then Dominator plays any other legal move. Also, if Staller decides to pass some additional moves (she still has $p-q$ pass moves at her disposal), this only helps Dominator to perform his strategy. This strategy thus guarantees that, no matter who was the last to play in the imagined game, the real game will be finished in no more than $2q$ additional moves, so that the total number of moves played in the real game is at most $(s-q)+ 2q = s+ q \le s+p$. Now, Dominator is playing optimally on the imagined game and Staller maybe not, hence $\gamma_g(G) \ge s$. On the other hand, Staller is playing optimally on the real game and Dominator maybe not, hence  ${\gamma}_g^{St,p}(G) \le s + q$. Therefore, ${\gamma}_g^{St,p}(G) \le s+q \le \gamma_g(G) + q \le \gamma_g(G) + p$.

We have thus proved the lemma under the assumption that every move of Staller from the real game is a legal move in the imagined game. Assume now that at some stage of the game the move $s_i$ of Staller in the real game is not legal in the imagined game and let $s^{(1)}, \ldots, s^{(k)}$, $k\in [q]$, be the vertices among the vertices $s^{(1)}, \ldots, s^{(q)}$ that were played in the imagined game until this moment. Let $X$ be the set of vertices dominated in the real game prior to Staller playing the vertex $s_i$. Since $s^{(1)}, \ldots, s^{(k)}$ are the only vertices played in the imagined game that were not played in the real game,
$$N[s_i]\setminus X\subseteq \cup_{j=1}^k N[s^{(j)}]\,.$$
Assume first that $N[s^{(\ell)}] \setminus (N[s_i] \cup X)\neq \emptyset$ for some $\ell\in [k]$. In this case Dominator replies to the move $s_i$ of Staller in the real game by playing the vertex $s^{(\ell)}$ (in the real game only).

Hence, at this stage, two additional moves ($s_i$ and $s^{(\ell)}$) were played in the real game, and then both games continue. Assuming that all the other moves of Staller from the real game are legal in the imagined game, the number of moves played in the real game at the time the imagined game ends is  $s - q + 2$. The strategy of Dominator to play in the real game the vertices  $s^{(j)}$, $j\in [k]\setminus \{\ell\}$, guarantees that the real game will finish in at most $(s-q+2) + 2(q-1) = s + q \le s + p$ moves.

The other case that can happen is that $\cup_{j=1}^k N[s^{(j)}] \subseteq N[s_i] \cup X$. In that case, however, already the move $s_i$ of Staller takes care of at least one of the vertices $s^{(j)}$ that are the target of Dominator after the imagined game finishes, so Dominator can play an arbitrary legal move (if there is such) in the real game and the conclusion from the end of the last paragraph holds.

Repeating the arguments from the above two paragraphs for possible later moves of Staller that are not legal in the imagined game completes the argument.
\qed

Lemma~\ref{lem:pass-game} has been proved for the case $p=1$ in~\cite[Lemma 2]{BKR2010}.

\section{Upper bound on the game domination number using edge cuts}
\label{sec:edge-cuts}

Let ${\cal C}(G)$ denote the set of minimal edge cuts of a graph $G$. Then our main result reads as follows.

\begin{theorem}
\label{thm:main1}
If $G$ is a connected graph, then
$$\gamma_g(G) \le \min_{C\in {\cal C}(G)}\{\gamma_g(G\setminus C) + 2|C|\}\,.$$
\end{theorem}

\proof
Let $C\in {\cal C}(G)$. Since $G$ is connected and $C$ is a minimal edge cut, the graph $G\setminus C$ consists of two components. Denote them with $G_1$ and $G_2$, so that $G\setminus C = G_1\cup G_2$. We thus need to prove that
\begin{equation}
\label{eq:key}
\gamma_g(G) \le \gamma_g(G_1\cup G_2) + 2|C|\,.
\end{equation}
Let a D-game be played on the graph $G$. In addition to that, Dominator imagines that a Staller $|C|$-pass game is played on $G_1\cup G_2$ as follows. Dominator plays optimally in the imagined game (played on $G_1\cup G_2$) and copies each such move to the real game (played on $G$), provided that the move is legal on $G$. Moreover, Dominator copies the moves of Staller from the real game (where Staller plays optimally) to the imagined game, provided that such a move of Staller is a legal move in the imagined game played on $G_1\cup G_2$. Suppose first that at some stage of the imagined game, a move of Dominator is not legal in the real game. In that case Dominator selects an arbitrary legal move in the real game provided there is such a move available; otherwise the game is over. Suppose next that at some stage of the real game, a move of Staller is not legal in the imagined game. Then in the imagined game Dominator assumes that Staller has passed her move.

We first claim that the game played on $G_1\cup G_2$ is indeed a Staller $|C|$-pass, game, that is, Staller makes at most $|C|$ pass moves. Suppose that a move of Staller, say $s_i$, in the real game is not a legal move in the imagined game.   Recall that from the strategy of Dominator described above,
 such a move corresponds to a pass by Staller in the imagined game.  Without loss of generality we assume $s_i \in V(G_1)$.
Since $s_i$ is not a legal move in the imagined game, the only newly dominated vertices of $G$ when Staller plays $s_i$ belong to $V(G_2)$.  Select any one of these
newly dominated vertices, say $v_i$, and assign the edge $s_iv_i$ to the vertex $s_i$.  Note that $s_iv_i$ belongs to $C$, and it now does not lie in the residual graph of the real game.
In this way each time Staller selects a vertex in the real game that is not legal in the imagined game, a new edge from $C$ is assigned to the move. This means that there is an
injection from the set of pass moves of Staller (imagined by Dominator) into $C$. This proves the claim.

In summary, the real game played is a D-game in which Dominator is maybe not playing optimally, while the imagined game is a Staller $|C|$-pass game in which Staller might not play optimally. Let $r$ be the number of moves played in the real game, $s$ the number of moves played in the imagined game, and $p$ the number of pass moves by Staller imagined by Dominator in the imagined game. Then, by the above claim, $p\le |C|$. Moreover, we claim that the following inequalities hold true:
\begin{enumerate}
\item[(i)] $\gamma_g(G)\le r$,
\item[(ii)] $r \le s + p$, and
\item[(iii)] $s\le \gamma_g^{St,|C|}(G_1\cup G_2)$.
\end{enumerate}
To see that (i) holds, recall that Staller plays optimally on $G$ but Dominator maybe not. Since the game has finished in $r$ moves it follows that $\gamma_g(G) \le r$.

(ii) Each move in the imagined game corresponds to a unique vertex played in the real game.  In addition, each (imagined) pass of Staller in the game played on $G_1 \cup G_2$
 was the result of a move by Staller in the real game (on $G$) that was not a legal move in the imagined game (on $G_1\cup G_2$). Since there are no additional moves in the real game
 and the real game finishes no later than the imagined game, it follows that $r \le s + p$.

(iii) The inequality $s\le \gamma_g^{St,|C|}(G_1\cup G_2)$ follows in a way similar to that stated in (i) because in the imagined game Dominator is playing the Staller $|C|$-pass
game optimally.

Putting together the fact $p\le |C|$ and claims (i)-(iii) we get
\begin{eqnarray*}
\gamma_g(G)  \le  r &\le&  s + p \\
& \le & \gamma_g^{St,|C|}(G_1\cup G_2) + p \\
& \le & \gamma_g(G_1\cup G_2) + |C| + p \\
& \le & \gamma_g(G_1\cup G_2) + 2|C|\,,
\end{eqnarray*}
where the before-last inequality follows by Lemma~\ref{lem:pass-game}.
\qed

The following consequences that follow directly from Theorem~\ref{thm:main1} are worth stating explicitly.

\begin{corollary}
\label{cor:edge-connectivity}
If $G$ is a connected graph and $C$ a minimum edge cut of $G$, then
$$\gamma_g(G) \le \gamma_g(G\setminus C) + 2\kappa'(G)\,.$$
\end{corollary}

\begin{corollary}
\label{cor:bridge}
If $e$ is a bridge of a connected graph $G$, then
$$\gamma_g(G) \le \gamma_g(G - e) + 2\,.$$
\end{corollary}

In~\cite[Theorem 2.1]{BDKK-2014} it was proved that Corollary~\ref{cor:bridge} holds for an arbitrary edge $e$ of a graph $G$, not only for bridges. In the (complicated) construction from ~\cite[Proposition 2.4]{BDKK-2014} it was also demonstrated that the equality can be achieved. The corresponding edge from it is not a bridge however.

\section{On game domination number and bridge removal}
\label{sec:bridge-removal}

In this section we show that the equality can be achieved in Corollary~\ref{cor:bridge}.  To accomplish this we introduce a  new concept.

We say that $G$ is a \emph{double-Staller graph} if it has the following property.  Suppose the Dominator $1$-pass game is played on $G$, that is, Dominator is allowed to pass at most once. Then Staller has a strategy that, for any strategy of Dominator, forces the game to last at least $\gamma_g(G)+p$ moves, where $p=0$ if Dominator does not pass and $1$ if he does. At  first sight this definition seems counterintuitive because why would Dominator pass a move if it prolonged the game by one move? However, if the game is played, say, on a disconnected graph, then Dominator might be forced to skip a move on some component and we could be interested in the game restricted to this component. Anyhow, we will see later how double-Staller graphs will be applied.

Examples of double-Staller graphs can be constructed as follows. For $n \geq 2$, let $T_n$ be the graph obtained from $K_{1,n}$ by attaching four leaves to each of the leaves of $K_{1,n}$. It is easy to check that $\gamma_g(T_n)=2n-1$.  On the other hand, if Dominator decides at some point to pass a move, then having enough leaves at her disposal Staller can prolong the game by one move.

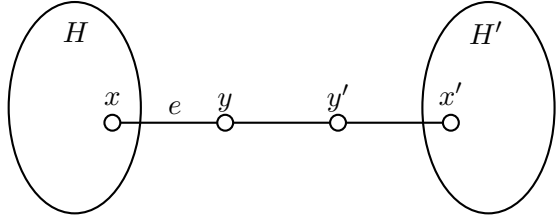
\begin{figure}[ht!]
\begin{center}
\begin{tikzpicture}[scale=1.0,style=thick,x=1cm,y=1cm]
\def\vr{3pt}
\coordinate(x) at (0,0);
\coordinate(y) at (1.5,0);
\coordinate(y') at (3,0);
\coordinate(x') at (4.5,0);
\draw (x) -- (y) -- (y') -- (x');
\draw (-0.5,0.2) ellipse (25pt and 40pt);
\draw (5,0.2) ellipse (25pt and 40pt);
\draw(x)[fill=white] circle(\vr);
\draw(y)[fill=white] circle(\vr);
\draw(y')[fill=white] circle(\vr);
\draw(x')[fill=white] circle(\vr);
\draw(x)++(0,0.1) node[above] {$x$};
\draw(y) node[above] {$y$};
\draw(y') node[above] {$y'$};
\draw(x')++(0,0.1) node[above] {$x'$};
\draw(0.6,0.2) node[right] {$e$};
\draw(-0.8,1.2) node[right] {$H$};
\draw(4.6,1.2) node[right] {$H'$};
\end{tikzpicture}
\end{center}
\caption{Graph $H_{x,3}$}
\label{fig:graphX}
\end{figure}

For an arbitrary graph $H$ and a vertex $x$ of $H$, let $H_{x,3}$ denote the graph of order $2n(H)+2$ constructed as follows.  Let $H'$ be an isomorphic copy of $H$ in which each vertex
$v\in V(H)$ is denoted by $v'$.  The graph $H_{x,3}$ is obtained from the disjoint union of $H$ and $H'$ by joining $x$ and $x'$ by a path of length 3.  The vertices on this path that are
adjacent to $x$ and $x'$ will be denoted by $y$ and $y'$, respectively.  Also set $e=xy$.  See Figure~\ref{fig:graphX}.

\begin{theorem} \label{thm:example}
Let $k$ be an even positive  integer and $H$ be a double-Staller graph.  If $H$ contains a vertex $x$ such that $\gamma_g(H|x)=\gamma_g(H)=k+1$ and
$\gamma_g'(H|x)=\gamma_g'(H)=k$, then $\gamma_g(H_{x,3})=2k+3$ and $\gamma_g(H_{x,3}-e)=2k+1$.
\end{theorem}

\proof
For ease of notation we let $G=H_{x,3}$ throughout the proof. We first prove that $\gamma_g(G)=2k+3$.

\medskip \noindent {\bf Claim A.} $\gamma_g(G)\leq 2k+3$. \\
Consider the following strategy of Dominator.  The first move of Dominator will be $d_1=y'$.
The general strategy of Dominator after playing $y'$ is to always respond to a move of Staller in $H$ (or in $H'$) so that the game restricted to $H$ (or restricted to $H'$) is the D-game or S-game without passes. We now distinguish two cases.

\medskip \noindent {\bf Case A.1.} $s_1=y$. \\
This move of Staller results in the residual graph being the disjoint union of $H|x$ with $H'|x'$.  In this case Dominator plays an optimal move $d_2$ in $H|x$ and then proceeds with his general strategy of ``following'' Staller in $H|x$ and $H'|x'$.  If the game on $H|x$ ends before the game on $H'|x'$, then since $\gamma_g(H|x)=\gamma_g(H)=k+1$ and $k+1$ is odd, Dominator made the last move on $H|x$; in particular, $k+1$ moves were made on $H|x$.  If any moves were made on $H'|x'$ before the game ended on $H|x$, then the last move made there was by Dominator.  Hence, after the game has ended on $H|x$ either the S-game starts on $H'|x'$ or the S-game continues there.  Since $\gamma_g'(H'|x')=k$,  the total number of moves played on $H'|x'$ is $k$.  Therefore, in this case the number of moves played was $2+(k+1)+k=2k+3$.

\medskip \noindent {\bf Case A.2.} $s_1 \ne y$. \\
Assume first that Staller never plays the vertex $y$.   Since $\gamma_g'(H)=\gamma_g'(H')=\gamma_g'(H'|x')$, we may assume without
loss of generality that $s_1 \in V(H)$.  Because of the general strategy of Dominator described above and the fact that $k$ is even, Staller will also be forced to be the first to play on $H'|x'$.  The result is that the S-game is played on both $H$ and $H'|x'$, which in turn implies only $2k+1$ moves are made.  Hence also in this case Dominator has achieved his goal of finishing the game in no more than $2k+3$ moves.

In the rest of this case we thus assume that at some point in the remainder of the game Staller plays $y$. Let $s_1\in V(H)$. We now distinguish the strategy of Dominator depending on when Staller plays $y$.

Suppose first that Staller plays $y$ before $H$ is dominated and before any vertex from $H'$ has been played. Then Dominator replies with an optimal move in $H$ and after that follows his general strategy. By Theorem~\ref{thm:diff-at-most-1}, the game on $H$ will last at most $k+1$ moves. If it lasts $k+1$ moves, then Staller is forced to play first on $H'$, so that $k$ moves will be played on $H'$. On the other hand, if the game lasted $k$ moves on $H$, then using Theorem~\ref{thm:diff-at-most-1} again, the game will last at most $k+1$ moves on $H'$. In any case, at most $2k+3$ moves will be made.

Note that after $H$ is dominated, the vertex $y$ is not a legal move.

Next, suppose that Staller plays $y$ after some vertices of $H$ and some vertices of $H'$ have been played, but neither $H$ nor $H'$ is dominated. Then Dominator plays his next move in $H$ and we can argue analogously as in the subcase when no vertex of $H'$ has yet been played that at most $2k+3$ vertices will be played.

In the final subcase assume that Staller plays $y$ after $H'$ is dominated and some vertices of $H$ have been played. Then Dominator replies with a move in $H$. Since $k$ vertices were played in $H'$, Theorem~\ref{thm:diff-at-most-1} again ensures at most $k+1$ moves will be played on $H$.

The above subcases cover all the possibilities when Staller plays her first move in $H$. If, however, she first plays in $H'$, then, since in all the cases considered she plays $y$ at some point, at that time the situation becomes symmetric (because both $y'$ and $y$ were played), hence the cases when she first plays in $H'$ can be treated along the same lines as those when she played first in $H$.

This proves Claim A.

\medskip \noindent {\bf Claim B.} $\gamma_g(G)\geq 2k+3$. \\
We next show that Staller has a strategy that ensures at least $2k+3$ moves are made on $G$.

\medskip \noindent {\bf Case B.1.} $d_1=y'$. \\
In this case Staller replies with $s_1=y$.
In the rest of the game the strategy of Staller is to respond to a move of Dominator in $H|x$ (or in $H'|x'$) with an optimal move for her in $H|x$ (or in $H'|x'$).

In the first subcase assume without loss of generality that $H$ was dominated before a single vertex from $V(H')$ was played.  Then, $k+1$ moves were played on $H$ and Dominator
made the last move there.  After this the S-game is played on $H'|x'$, which requires $k$ additional moves.  Hence, $2+(k+1)+k=2k+3$ moves were played as claimed.

In the second subcase $H$ was dominated before $H'$, but some vertices from $H'$ had already been played.  Because of Staller's strategy to follow Dominator, the last of these
moves in $H'$ was made by Staller.   As in the first subcase, $k+1$ vertices were played in $H$, and it is now Staller's move, which will be in $H'$.
Suppose $2\ell$ moves were made in $H'$ before $H$ was dominated.  Let $S$ be the set of vertices from $H'$ that are dominated by these $2\ell$
vertices. It is now Staller's turn to play on the partially dominated graph $H'|(S\cup \{x'\})$.  Then, by Corollary~\ref{cor:prevertex} and by the fact that $H'$ is a double-Staller graph,
$$\gamma_g'(H'|(S \cup \{x'\})) \ge \gamma_g'(H'|S)-2 \geq ((k+2)-2\ell)-2\,,$$
hence $\gamma_g'(H'|(S \cup \{x'\})) \ge k-2\ell.$  This in turns means that at least $k$ moves will be played on $H'$, and hence $2k+3$ vertices in total were played.

The situation when $d_1=y$ is the same.

\medskip \noindent {\bf Case B.2.} $d_1=x$. \\
Now the strategy of Staller is to play optimally in $H$ so long as Dominator continues to play vertices from $H$.  If both players
continue to make moves in $H$ until $H$ is dominated, then $k+1$ moves were made and Dominator made the last move.  Staller then plays the vertex $y$.  By the Continuation Principle the next
move of Dominator will be on $H'$ and the strategy of Staller to play on $H'$ guarantees that at least $k+1$ moves were played on $H'$  when the game is finished.  Thus,
a total of $(k+1)+1+(k+1)=2k+3$ moves were made.

Suppose next that at some
point when Staller made a move on $H$, Dominator played a vertex in $V(H') \cup \{y'\}$.  By the Continuation Principle, Dominator will not play $y'$.  Staller's next move will
be on $H'$.  After this move, no matter how Dominator switches between playing a vertex from $H$ or a vertex from $H'$, Staller follows her optimal strategy on that subgraph.  We may assume
without loss of generality that $H$ is dominated first.  In this case, $k+1$ moves were made on $H$.  Since Dominator was the last to play on $H$, Staller now has two consecutive moves on $H'$.  Since $H'$ is a double-Staller graph, at least $k+2$ moves will be made on $H'$, and consequently at least $2k+3$ moves in total.

The situation when $d_1=x'$ is the same.

\medskip \noindent {\bf Case B.3.} $d_1\in V(H)\setminus \{x\}$. \\
Staller again employs the strategy of playing optimally in $H$ as long as Dominator plays there.  If all of $H$ is dominated before any vertex
outside of $H$ is played, then $k+1$ vertices were played in $H$, and it is Staller's turn. She plays $y$ as her next move.  By the Continuation Principle, Dominator replies with a move in $H'$, which
yields at least $(k+1)+1+(k+1)=2k+3$ moves when the game ends.  Otherwise (before all of $H$ is dominated), Dominator plays a vertex $u$ not in $H$.  Note that $u$ is the first vertex
outside of $H$ that was played by either of the players.  Now we distinguish three subcases.

\medskip \noindent {\bf Case B.3.1.} $u=y$. \\
Staller replies by playing $y'$.  Now, by arguing as in the second subcase of Case B.1, it will follow that at least $2k+3$ vertices will be played when the game ends.

\medskip \noindent {\bf Case B.3.2.} $u=y'$. \\
Staller replies with an optimal move in $H$, and after that Staller follows Dominator in $H$ or $H'$.  (Note that by the Continuation Principle, we may assume that
Dominator will never play $y$.) Suppose that $H$ is dominated before $H'$.   Since $H$ is a double-Staller graph, $k+2$ is even, and Staller played two consecutive moves in $H$, it follows
that Dominator was the last to play in $H$. Suppose first that when Staller plays the next move on $H'$, it is in fact the first vertex played there.  In this case the S-game is played on
$H'|x'$, and hence exactly $k$ moves will be made there.  Thus, $(k+2)+1+k=2k+3$ are played in the game.  Suppose next that when Staller makes her first move just after $H$ was dominated
some vertices had already been played from $H'$.  But then we can argue again just as in the second subcase of Case B.1 that $k$ moves will be played on $H'$, which again shows that
$2k+3$ moves are forced by Staller.

Assume now that $H'$ was dominated before $H$.  In this case $k+1$ moves were played in $H'$, and Dominator made the final move there.  Hence, the rest of the game is the S-game played on
$H|S$, where $S$ is the set of vertices already dominated in $H$.  If the D-game were played on $H|S$, then $k+2$ moves would be played on $H$ since $H$ is a double-Staller graph.  By
Theorem~\ref{thm:diff-at-most-1} it follows that at least $k+1$ vertices will be selected in $H$.  Therefore, at least $2k+3$ vertices will be played in total.

\medskip \noindent {\bf Case B.3.3.} $u\in V(H')$. \\
If Dominator plays all the time on $V(H) \cup V(H')$, then on one of $H$ or $H'$ we will have $k+1$ moves, while on the other Staller will have two
consecutive moves.  Since $H$ and $H'$ are double-Staller graphs, the result is that at least $2k+3$ moves are guaranteed.

Assume next that Dominator plays $y$ before either $H$ or $H'$ is dominated.  (The possibility that Dominator plays $y'$ is treated in the same manner.) Staller now makes an optimal
move in $H'$, and then follows Dominator in $H$ and $H'$.  If $H'$ is dominated before $H$, then $k+2$ moves were made there and Dominator made the last one.  Once again, as in the
last subcase of Case B.1, we see that at least $k$ vertices were selected in $H$.  If $H$ is dominated before $H'$, then $k+1$ moves were made on $H$ and Dominator was the last to play
on $H$.  Staller has the next move on $H'$, and by the same argument as in the last paragraph of Case B.3.2 we see that at least $k+1$ moves will be made in $H'$, which once again implies
that at least $2k+3$ moves are made in total.

Finally, assume that one of $H$ or $H'$ was dominated before $y$ was played (if $y$ was played at all).  Suppose that it is $H$ that was dominated before $H'$.   Note that the
  last move on $H$ was made by Dominator.  If $x'$ was not already played, then $y'$ is a legal move and Staller plays it.  In this way Staller has forced Dominator to make the next move in $H'$, which implies that $k+1$ vertices will be played there.  This gives a total of at least $2k+3$ vertices played when the game ends.  On the other hand if $x'$ was already played, then Staller replies with an optimal move in $H'$.  If Dominator does not play the vertex $y$ in the remainder of the game, then since $H'$ is a double-Staller graph, we will again have at least $k+2$ moves in $H'$.
  On the other hand, if at some point Dominator plays $y$, then by the Continuation Principle there will still be at least $k+1$ vertices selected from $H'$.  Together with $y$ and the $k+1$
  vertices played in $H$ this shows that at least $2k+3$ moves were made when the game ends.  The other possibility, that $H'$ was dominated before $H$, can be treated in a similar way.

This proves Claim B.

\medskip To complete the proof we must show that $\gamma_g(G-e)=2k+1$.  Since $\gamma_g(G)=2k+3$ has already been proved, the result of \cite[Theorem 2.1]{BDKK-2014} implies that $\gamma_g(G-e)\geq 2k+1$. Hence, we only need to provide a strategy for Dominator that will guarantee no more than $2k+1$ moves will be made when the D-game is played on $G-e$. Dominator starts the game with $d_1=y'$ and in the rest of the game follows the moves of Staller on $H$ and $H'|x'$.  In this way, since $k$ is even, Staller will be forced to make the first move in both $H$ and $H'|x'$.  Since $\gamma_g'(H'|x')=k=\gamma_g(H)$, the game will end in precisely $2k+1$ moves.  \qed

Examples of graphs $H$ that give the claimed conclusions of Theorem~\ref{thm:example} include $C_6$, $C_{10}$ and $C_{14}$.  These were checked by computer.  Hence it seems reasonable
to conjecture that these conclusions hold for all graphs in the infinite family $C_{4k+2}$, $k \ge 1$.

\section{An application}
\label{sec:application}

Our original motivation for considering edge cuts in the frame of the domination game was \cite[Conjecture 1.1]{james-2019} asserting that if $G$ is traceable, then $\gamma_g(G) \le \lceil n(G)/2\rceil$.  Note that by the results of~\cite{kosmrlj-2017} this conjecture holds for paths and cycles.  As far as we know there is no other nontrivial class of traceable graphs for which this conjecture has been verified.  With the following theorem new infinite families can be constructed for which the conjecture holds.

\begin{theorem}
\label{thm:hamiltonian}
Let $C$ be a minimum edge cut of a connected graph $G$ and let $G_1$ and $G_2$ be the  components of $G\setminus C$. If $\gamma_g(G_2) \leq q\cdot n(G_2)$ and
$\gamma_g(G_1) \le \frac{n(G_1)}{2} - \frac{(2q-1)n(G_2)}{2} - 2\kappa'(G) - 2$, then $\gamma_g(G) \le \lceil n(G)/2\rceil$.
\end{theorem}

\proof
 From~\cite[Theorem 3.3]{dorbec-2015} we  know that $\gamma_g(G_1\cup G_2) \le \gamma_g(G_1) + \gamma_g(G_2) + 2$.  Using this fact and Corollary~\ref{cor:edge-connectivity} we can compute as follows:
\begin{eqnarray*}
\gamma_g(G) & \le & \gamma_g(G_1\cup G_2) + 2\kappa'(G) \\
& \le & (\gamma_g(G_1) + \gamma_g(G_2) + 2) + 2\kappa'(G) \\
& \le & \left(\frac{n(G_1)}{2} - \frac{(2q-1)n(G_2)}{2} - 2\kappa'(G) - 2\right) + q\cdot n(G_2) + 2 + 2\kappa'(G) \\
& = & \frac{n(G_1)}{2} + \frac{n(G_2)}{2} \\
& \le & \left\lceil\frac{n(G)}{2}\right\rceil\,,
\end{eqnarray*}
which we wanted to show.
\qed

The best published upper bound for the game domination number for a graph of order $n$ with no isolated vertices is  $\gamma_g(G) \le 2n(G)/3$;
 see~\cite[Proposition~2]{bujtas-2015}.
Assuming that the graph $G_2$ in Theorem~\ref{thm:hamiltonian} has order at least two the condition on $\gamma_g(G_1)$ reduces to
\begin{equation} \label{eqn:2/3}
\gamma_g(G_1) \le \frac{n(G_1)}{2} - \frac{n(G_2)}{6} - 2\kappa'(G) - 2.
\end{equation}
In this respect we mention that an improvement of the upper bound from $\frac{2n(G)}{3}$
to $\frac{5n(G)}{8}$ has been announced (see~\cite[p.4]{bujtas-2018+}), in which case the condition on $\gamma_g(G_1)$ can be further weakened.
Further improvements are possible if $G_2$ has no pendant vertices~\cite{heki-2016}, if the minimum degree of $G_2$ is at least $3$~\cite{bujtas-2015}, and also if $G_2$ is a so-called weakly $S(K_{1,3})$-free forest~\cite{schmidt-2016}.

To see how Theorem~\ref{thm:hamiltonian} can be applied consider the following infinite family of traceable graphs.  Let $H$ be any traceable graph of order at most $3k-42$ where
$k \ge 15$.  Let $a$ and $b$ be the ends of a hamiltonian path in $H$.  Select two vertices $u$ and $v$ of a complete graph of order $k$. Let $K_k(H)$ be the graph constructed from the disjoint union of $K_k$ and $H$ by adding the edges $au$ and $bv$.  Note that the set $\{au,bv\}$ is a minimum edge cut of $K_k(H)$ and that $K_k(H)$ is traceable. By setting $G_1=K_k$
and $G_2=H$ it is easy to check that \eqref{eqn:2/3} is fulfilled, and hence Theorem~\ref{thm:hamiltonian} applies.

\section*{Acknowledgements}

The authors are very grateful for the anonymous referee who helped us to make precise the notions of a Staller $p$-pass game and a double-Staller graph.
The financial support from the Slovenian Research Agency is acknowledged (research core funding No.\ P1-0297, projects N1-0043, J1-9109, and the bilateral grant BI-US/18-19-011).
Both authors acknowledge the support of the Wylie Enrichment Fund of Furman University.

\end{document}